\documentclass[11pt, english]{amsart}
\usepackage{babel,amssymb,amsmath,latexsym,amsthm,amscd,eucal,enumerate,verbatim}  
\usepackage[dvips]{graphicx}
\usepackage[all]{xy}
\usepackage[left=2.5cm,top=2.5cm, bottom=2.5cm,right=2.5cm]{geometry}
\usepackage{setspace}

\newtheorem{thm}{Theorem}

\newtheorem{lem}[thm]{Lemma}

\newcommand\bp{\noindent{\it Proof.}\ }

\title{Positive maps which map the set of rank k projections onto itself}

\author{Erling St{\o}rmer}

\date{18-4-2016}
\begin{document}
\maketitle

\begin{abstract}

{Extending Wigner\rq s theorem we give a characterization of positive maps of $B(H)$ into itself which map the set of rank k projections onto itself}
\end{abstract}

  One form of the celebrated Wigner\rq s theorem \cite{4W} is that if $ \phi$ is a linear map of the bounded operators $B(H)$ on a Hilbert space $H$ into itself with the property that it maps the set of rank 1 projections bijectively onto itself, then $\phi$ is of the form

\begin{equation}                      \phi(a) = UaU^* \quad   \text{or} \quad    \phi(a) = Ua^t U^*, \tag{*} \end{equation}

\medskip\noindent
where $a^t$ is the transpose of $a$ with respect to a fixed orthonormal basis for  $H$, and $U$ is a unitary operator.  In the paper \cite{Sa}  Sarbicki, Chruscinski and Mozrzymas generalized this to the case when $H$ is of finite dimension n with n a prime number, and the set of rank 1 projections is replaced by rank k projections, where k is a natural number strictly smaller than n. They gave a counter example to the conclusion (*) when n is not a prime.  In that case $\phi$ is no longer a positive map. In the present note we make the extra assumption that $\phi$ is a positive unital map.  Then for any Hilbert space we obtain the conclusion (*). Closely related results have been obtained by Molnar \cite{M}.
  
  Recall that an atomic masa in $B(H)$ is a maximal abelian subalgebra $A$ generated by the rank 1 projections  corresponding to the vectors in an orthonormal basis for $H$.  Thus if $H$ is finite dimensional each maximal abelian subalgebra is atomic. We start with a lemma.  See also \cite{M}, Lemma 2.1.5.
 
 \begin{lem}\label {lem}
 Let $p\in B(H)$ be a rank 1 projection and $A$ an atomic masa in $B(H)$ containing $p$.  Let $k$ be a natural number, $k < dim H$. Then there exist $k+1$ projections $P_1,...,P_{k+1}$ in $A$ such that
 $$
 p = \frac{1}{k} \sum _{j=2}^{k+1} P_j  - \frac{k-1}{k} P_1.$$

 \end{lem}

 \bp
 Let $p_1=p ,p_2,...,p_{k+1}$ be mutually orthogonal rank 1 projections in $A$. Let
 $$
 P_j = \sum_ {i=1,i\neq j}^{k+1} p_i,  j = 1,...k+1.
 $$
 Then $P_j$ is a projection of rank k, and $p_i\leq P_j$ for all $j\neq i$,  so $p_i \leq  P_j$ for $k$ of the projections $  P_j$. It is therefore an easy computation to show the above formula.  The proof is complete.
 
 \begin{thm}\label{thm 2}

 Let $\phi$ be a positive unital map of $B(H) $ into itself such that $\phi$ maps the set of  projections of rank $k$ in $B(H), k<dim H$, onto itself.   Then $\phi$ is of the form (*).

 \end {thm}
 \bp
 Since each projection of rank k in $B(H)$ is in the image under $ \phi$  of a rank k projection, it follows from Lemma 1 that the rank 1 projections are in image of $\phi$, hence each finite rank operator is in the image of the finite rank operators.  By continuity of $\phi$ it follows that $\phi$ when restricted to the compact operators $C(H)$, maps $C(H)$ onto a norm dense subset of itself.

The definite set $D$ of $\phi$ is the set of self-adjoint operators $a$ such that $\phi(a^2) = \phi(a)^2$.  Let $Q$ be a projection of rank k; then $P = \phi(Q)$ is a projection of rank k, hence
$$
\phi(Q^2) = \phi(Q) = P = P^2 = \phi(Q)^2,
$$
so that $Q \in D$. By  \cite{St}, Proposition 2.1.7, $D$ is a norm closed Jordan subalgebra of $B(H)$, so by the same argument as above $D\cap C(H)  = C(H)_{sa} $ , the self-adjoint operators in $C(H)$. Furthermore, the restriction of $\phi$ to $D$ is a Jordan homomorphism.  Since $C(H)$ is irreducible, by \cite{S}, Corollary 3.4, $\phi$ is either a homomorphism or an anti-homomorphism on $C(H)$.  But $C(H)$ is a simple C*-algebra, so $\phi$  is either an automorphism or an anti-automorphism of $C(H)$.
 
 Let now $\omega_x$ be a vector state on $B(H)$.  Then if $p$ is the rank 1 projection onto the 1-dimensional subspace of $H$ generated by $x$ , then for $a\in B(H)$,
  $$
  \omega_x (a) = (ax,x) = Tr(pap).
  $$
  Since $\phi$ is a Jordan automorphism of $C(H)$ there is a unit vector $y$ such that if $q$ is the rank 1 projection onto the subspace spanned by $y$, then $\phi(q) = p$.
  Thus for $a\in B(H)$, since $q\in D$,  we have, using  \cite{St}, Proposition 2.1.7, 
  $$
  \omega_x (\phi(a)) = Tr(p \phi(a)p) = Tr(\phi(q)\phi(a)\phi(q)) = Tr(\phi(qaq)) = Tr(\phi(\omega_y (a))q) = \omega_y(a)Tr(p) = \omega_y(a).
  $$
We have thus shown that each vector state composed with $\phi$ is a vector state.  Hence by \cite{St} ,Theorem 3.3.2 $\phi$ is  of the desired form (*).  The proof is complete.

Department of Mathematics,  University of Oslo, 0316 Oslo, Norway.

e-mail  erlings@math.uio.no

\end{document}